\documentclass[11pt]{article}
\usepackage{latexsym,amssymb,amsmath,a4wide,theorem}

\newtheorem{theorem}{Theorem}[section]
\newtheorem{lemma}{Lemma}[section]
\newtheorem{proposition}[theorem]{Proposition}

\newcommand{\f}{\frac}
\newcommand{\lf}{\left}
\newcommand{\rg}{\right}
\newtheorem{rmk}{Remark}[section]

\title{Stability of degenerate stationary solution to the outflow problem for full Navier-Stokes
equations}

\author{Yazhou Chen ,
\thanks{Department of Mathematics, School of Science, Beijing University of Chemical Technology, Beijing 100029,
P R China (chenyz@mail.buct.edu.cn).}\\
{\small } \and Hakho Hong ,
\thanks{Institute of Mathematics, State Academy of Sciences, Pyongyang, D P R Korea(hhong@amss.ac.cn).}\\
{\small } \and Xiaoding Shi
\thanks{Corresponding author. Department of Mathematics, School of Science, Beijing University of Chemical Technology, Beijing 100029,
P R China (shixd@mail.buct.edu.cn).}\\
}
\date{}
\begin{document}
\maketitle
\begin{abstract}
This paper is concerned with the large-time behavior of solutions to
the outflow problem of full compressible Navier-Stokes equations in
the half line. This is one of the series of papers by the authors on
the stability of nonlinear waves to the outflow problem. We show the
time asymptotic stability of degenerate (transonic) stationary
solution for the general gas including ideal polytropic gas.
The key point of the proof is based the new property of the degenerate stationary solution and the delicate energy estimates.

\vspace{.20cm}\noindent\textbf{MSC 2010:} 35Q30, 35B35, 35L65,
76D33, 74J40.

 \vspace{.20cm} \noindent\textbf{Keywords:} compressible
Navier-Stokes equations, outflow problem, stationary solution,
 stability.
\end{abstract}

\section{Introduction}
 \setcounter{equation}{0}
The one-dimensional compressible Navier-Stokes equations are described in Eulerian coordinate by the system
\begin{equation}\label{1.1}\begin{cases}
\rho_{t}+(\rho u)_{x}= 0, \qquad x>0,\,\,\, t>0,\\
(\rho u)_{t}+(\rho u^2+p)_{x}
=\mu u_{xx},\\
[\rho(e+\frac{u^{2}}{2})]_t+[\rho u(e+\frac{u^{2}}{2}) +pu]_{x}
=\kappa\theta_{xx}+\mu(uu_{x})_{x},
\end{cases}
\end{equation}
where $u(x,t)$ is the velocity, $\rho(x,t)>0$ is the density,
$\theta(x,t)>0$ is the absolute temperature, $p=p(\rho,\theta)$ is
the pressure and $e=e(\rho,\theta)$ is the internal energy, while
$\mu$ and $\kappa$ denote the viscosity and the heat-conductivity
respectively.
Throughout this paper, the pressure $p(\rho,\theta)$ and the
internal energy $e(\rho,\theta)$ are assumed to satisfy
\begin{equation}\label{1.8-3}
p_{\rho}(\rho,\theta)>0,\hspace{0.6cm} e_{\theta}(\rho,\theta)>0.
\end{equation}
We consider the initial-boundary value (IBV) problem to the system
\eqref{1.1} on $[0, \infty)$ supplemented with the initial data
\begin{equation}\label{1.2}\begin{cases}
(\rho,u,\theta)(x,0)=(\rho_0,u_0,
 \theta_0)(x),\quad x>0,\\
 \displaystyle\lim_{x\rightarrow\infty}(\rho,u,\theta)(x,t)=(\rho_+,u_+,\theta_+),\quad t\geq 0,\end{cases}
\end{equation}
 and the boundary condition
\begin{equation}\label{1.3}
 u(0,t)=u_-<0,\qquad \theta(0,t)=\theta_-,
\end{equation} where $\theta_\pm>0,\rho_+>0,u_\pm
$ are prescribed constants.
\noindent\begin{rmk}
According to the sign of the velocity $u_-(=0,>0,<0)$ on the boundary $x=0$,
the following three type of problems are proposed \cite{M 99}: the impermeable wall problem, the inflow problem and the outflow problem.
It should be noted that for the inflow problem, the density $\rho_->0$ must
be given on the boundary by the well-posedness theory of the
hyperbolic equation $(1.1)_1$.
\end{rmk}

For the Cauchy problem of \eqref{1.1}, it is well known that the large time behavior of the solutions are described by the corresponding Riemann solutions to the hyperbolic part of the system \eqref{1.1} which consist of three basic wave patterns (shock wave, rarefaction wave and contact discontinuity) and  their superpositions in the increasing order of characteristic speed. But in the case of  the IBV problem of \eqref{1.1}, not only basic wave patterns but also a stationary solution, which is called the boundary layer solution (BL-solution for brevity), may appear due to
the boundary effect. For the IBV problem of isentropic Navier-Stokes
equations, Matsumura \cite{M 99} proposed a criterion on the
question when the BL-solution forms  and  a complete classification
about the precise description of the large time behaviors of
solutions. Since then, many results have been obtained for the
rigorous mathematical justification of this classification. We refer
to Matsumura-Mei \cite{MM 99}, Matsumura-Nishihara \cite{MN 00} etc. for the impermeable wall problem, to Kawashima-Nishibata-Zhu \cite{KNZ
03}, Nakamura-Nishibata-Yuge \cite{NNY 07}, Kawashima-Zhu \cite{KZ 08,KZ 09}, Huang-Qin \cite{HQ 09} etc. for the outflow problem, and to
Matsumura-Nishihara \cite{MN 01}, Huang-Matsumura-Shi \cite{HMS 03}, Shi \cite{S 03}, Fan-Liu-Wang-Zhao \cite{FLWZ 14} etc. for the inflow problem and so on. Further, for non-isentropic Navier-Stokes equations (i.e. problem \eqref{1.1}), we refer to Huang-Li-Shi \cite{HLS 10} etc.
for the impermeable wall problem, to Kawashima-Nakamura-Nishibata-Zhu \cite{KNNZ 10}, Qin \cite{Q 11}, Wan-Wang-Zhao \cite{WWZ 15-1, WWZ
15-2}, Chen-Hong-Shi \cite{CHS 17} etc. for the outflow problems and to Qin-Wang \cite{QW 09, QW 11}, Nakamura-Nisibata \cite{NN 11}, Zheng-Zhang-Zhao \cite{ZZZ 11}, Hong-Wang \cite{HW 16, HW 16-1} etc. for the inflow problems.  More works, please refer to the references therein.

The following focus on the more details for the stability and
convergence rate toward stationary solutions to the outflow problem
concerned with this paper.
For the isentropic Navier-Stokes equations, Kawashima-Nishibata-Zhu \cite{KNZ 03} first proved asymptotic stability of stationary
solutions under small $H^1-$initial perturbation. Kawashima-Zhu
\cite{KZ 08} improved the result in \cite{KNZ 03} to the combination
of stationary solution and rarefaction wave,  and Huang-Qin \cite{HQ
09} improved the results in \cite{KNZ 03, KZ 08} to large initial
perturbation. For this stability theorem, a convergence rate was
obtained by Nakamura-Nishibata-Yuge \cite{NNY 07} by assuming that
an initial perturbation belongs to the weighted Sobolev space.
Nakamura-Ueda-Kawashima \cite{NUK 09} gave a refinement of the
stability result established in \cite{NNY 07} for a degenerate
stationary solution. Precisely they obtained a convergence rate
under a more moderate assumption rather than the result in \cite{NNY
07}. The generalization of these one-dimensional problem to a
multi-dimensional half space was also studied, i.e. Kagei-Kawashima
\cite{KK 06} proved the asymptotic stability of a planar stationary
solution, and Nakamura-Nishibata \cite{NN 09} obtained the
convergence rate.
For full Navier-Stokes equations \eqref{1.1} of ideal polytropic
gas, Kawashima-Nakamura-Nishibata-Zhu \cite{KNNZ 10} first proved
the existence, the asymptotic stability and the convergence rate of
the stationary solution. Qin \cite{Q 11} proved that the
non-degenerate (supersonic or subsonic) stationary solution is
asymptotically stable under ¡Èpartially¡É large initial perturbation
with the technical condition that the adiabatic exponent $\gamma$ is
close to 1. Just recently, Wan-Wang-Zou \cite{WWZ 15-1}  improved
the result in \cite{Q 11} to the large initial perturbation without
any restriction on the adiabatic exponent $\gamma$. Also,
Wan-Wang-Zhao \cite{WWZ 15-2} studied the case when the
corresponding time-asymptotic state is a rarefaction wave or its
superposition with a non-degenerate stationary solution under large
initial perturbation. Just recently, Chen-Hong-Shi \cite{CHS 17} studied
the stability and convergence rate of a non-degenerate stationary
solution to outflow problem in the case of general gas \eqref{1.8-3}
including ideal polytropic gas, that is
\begin{equation}\label{perfect}
p=R\rho\theta=A\rho^{\gamma}e^{\f{\gamma-1}{R}s}\quad
\text{and}\quad e=\frac{R}{\gamma-1}\theta,
\end{equation} where $\gamma>1$ is the adiabatic exponent and $A, R$
are both positive constants.

\vspace{0.2cm}Although considerable progress has been obtained for
the stability of stationary solutions to the outflow problem of full
Navier-Stokes system \eqref{1.1}, however most of the results are
obtained only in the case of the non-degenerate stationary solution
except for Kawashima-Nakamura-Nishibata-Zhu \cite{KNNZ 10} where they proved the stability of
degenerate stationary solution for the ideal polytropic gas
\eqref{perfect}, but for the general gas, it is not trivial and the implicit realtions of various thermodynamical variables lead to many complicated terms in the course of establishing the energy estimate.

In this paper, we are interested in  the stability of the
degenerate (transonic) stationary solution of the outflow problem of \eqref{1.1}-\eqref{1.3} for the general gases satisfying
\eqref{1.8-3} (see Theorem \ref{theo 1.1}).
The main difficulty  for this problem is
that the degenerate stationary solution has only the algebraic decay property, not the exponential decay property, therefore  the effective methods which used to deal with the non-degenerate (supersonic and subsonic) case (see \cite{Q 11, WWZ 15-1,WWZ 15-2, CHS 17} etc.) are no longer applicable. The key point of the proofs in this paper is to derive the estimates (2.4) in Lemma 2.1, and this is mainly due to our key Proposition 1.1 below which enables us to control the lower order term $I_i(i=1,2,3,4,5)$ in Lemma 2.1.


\vspace{0.2cm} Now we will give some results on the stationary solutions to the outflow problem of full Navier-Stokes system and state our main result. Setting $v=\rho^{-1}$, it is well-known that by using
any given two of the five thermodynamical variables ($v, p, e,$ $
\theta$ and $s$), the remaining three variables are smooth functions
of them, where $s$ is the entropy of the gas. The second law of
thermodynamics $\theta ds=de+pdv$ asserts that, if we choose
$(v,\theta)$ or $(v,s)$ as independent variables and write
$(p,e,s)=(p(v,\theta),e(v,\theta),s(v,\theta))$ or
$(p,e,\theta)=(\widetilde{p}(v,s),\widetilde{e}(v,s),\widetilde{\theta}(v,s))$,
respectively, then we can deduce that
\begin{equation}\label{1.6-3}\begin{cases}
  s_{v}(v,\theta)=p_{\theta}(v,\theta), \\
s_{\theta}(v,\theta)=\frac{e_{\theta}(v,\theta)}{\theta},\\
e_{v}(v,\theta)=\theta p_{\theta}(v,\theta)-p(v,\theta),\\
\end{cases}
\end{equation}
or
\begin{equation}\label{1.7-3}\begin{cases}
  \widetilde{e}_{v}(v,s)=-\widetilde{p}(v,s),\qquad \qquad \qquad \widetilde{e}_{s}(v,s)=\theta,\\
\widetilde{p}_{v}(v,s)=p_{v}(v,\theta)-\frac{\theta
(p_{\theta}(v,\theta))^{2}}{e_{\theta}(v,\theta)},
\quad \widetilde{p}_{s}(v,s)=\frac{\theta p_{\theta}(v,\theta)}{e_{\theta}(v,\theta)}, \\
\widetilde{\theta}_{v}(v,s)=-\frac{\theta
p_{\theta}(v,\theta)}{e_{\theta}(v,\theta)}, \qquad \qquad \qquad
\widetilde{\theta}_{s}(v,s)=\frac {\theta}{e_{\theta}(v,\theta)}.
\end{cases}
\end{equation}
From \eqref{1.7-3} and \eqref{1.8-3} we have
\begin{equation}\label{1.9-3}
 \widetilde{p}_{v}(v,s)=p_{v}(v,\theta)-\frac{\theta
(p_{\theta}(v,\theta))^{2}}{e_{\theta}(v,\theta)}<0,
\end{equation}
\begin{equation}\label{1.10-3}\begin{cases}
 \widetilde{e}_{ss}(v,s)=\frac {\theta}{e_{\theta}(v,\theta)}>0,\quad
 \widetilde{e}_{vs}(v,s)=-\frac{\theta p_{\theta}(v,\theta)}{e_{\theta}(v,\theta)}, \\
 \widetilde{e}_{vv}(v,s)=-p_{v}(v,\theta)+\frac{\theta
(p_{\theta}(v,\theta))^{2}}{e_{\theta}(v,\theta)}>0,
\end{cases}\end{equation}
which means that $\widetilde{e}(v,s)$ is convex with respect to
$(v,s).$

\vspace{0.2cm}\textbf{Notation:} Throughout the rest of this paper,
$O(1), c$ or $C$ will be used to denote a generic positive constant
independent of  $x, t$ and $c_{i}(\cdot, \cdot)$ or $C_{i}(\cdot,
\cdot)(i\in Z_{+})$ stands for some generic constants depending only
on the quantities listed in the parentheses. As long as no confusion
arises, we denote the usual Sobolev space with norm
$\|\cdot\|_{H^{k}}$ by $H^k:=H^k(0,\infty)$ and
$\|\cdot\|_{H^{0}}=\|\cdot\|$ will be used to denote the usual
$L_2-$norm.

\vspace{0.2cm}  The stationary solution $(\hat{\rho}, \hat{u},
\hat{\theta})(x)$ of the outflow problem \eqref{1.1}-\eqref{1.3}
must satisfy the following system:
\begin{equation}\label{1.12-3}\begin{cases}
&(\hat{\rho} \hat{u})_{x}= 0, \qquad x>0,\\
&(\hat{\rho} \hat{u}^2+\hat{p})_{x}
=\mu \hat{u}_{xx},\\
&[\hat{\rho}\hat{u}(\hat{e}+\frac{\hat{u}^{2}}{2})
+\hat{p}\hat{u}]_{x}
=\kappa\hat{\theta}_{xx}+\mu(\hat{u}\hat{u}_{x})_{x},
\end{cases}
\end{equation}
\begin{equation}\label{1.12-31}
(\hat{u}, \hat{\theta})(0)=(u_-, \theta_-),\quad
\lim_{x\rightarrow\infty}(\hat{\rho}, \hat{u},
\hat{\theta})(x)=(\rho_+, u_+, \theta_+),
\end{equation}
where $\hat{p}=p(\hat{\rho},
\hat{\theta}),\,\,\,\hat{e}=e(\hat{\rho}, \hat{\theta}).$

Define the sound speed $c(v,s)$ and the Mach number $M(v,u,s)$ is
defined, respectively, by
$$c(v,s)=\sqrt{\frac{\partial p(\rho,s)}{\partial\rho}}=\sqrt{-v^2\widetilde{p}_v(v,s)},
\quad M(v,u,s)=\frac{|u|}{c(v,s)},$$ and set $M_+=\frac{|u_+|}{c(v_+,
s_+)},$ where $s=s(v,\theta)$ and $s_+=s(v_+,\theta_+).$

\vspace{0.2cm}The existence and the properties of the stationary solution
$(\hat{\rho}, \hat{u}, \hat{\theta})(x)$ satisfying \eqref{1.12-3}
and \eqref{1.12-31} are quoted in the following proposition which is proved by \cite{CHS 17}.
The proof is given in the Appendix for reader's convenience.

\begin{proposition} \label{theo 1.1-3}(Existence of stationary solution). Let
$\rho_+>0,\,\,u_-<0,\,\,\theta_\pm>0$. The necessary condition for
the existence of a solution to the system \eqref{1.12-3},
\eqref{1.12-31} is
\begin{equation}\label{2.3-3} \hat{\rho} \hat{u}= \rho_+ u_+=
\hat{\rho}(0) u_-,\,\,\forall x>0.
\end{equation}
So, if $u_+\geq 0$, there is no solution to the system
\eqref{1.12-3}, \eqref{1.12-31}.
If $u_+<0$ and \eqref{1.8-3} is hold, then we have the following
results.

1) For supersonic case $M_+>1$, there exists a positive constant
$\delta_0$  such that if $(u_-, \theta_-)\in
\mathcal{M}_{\delta_0}:=\{(u,\theta)\in R_-\times R_+\mid
|(u-u_+,\theta-\theta_+)|\leq \delta_0\}$, then there exists a
unique solution  $(\hat{\rho}, \hat{u}, \hat{\theta})(x)$ satisfying
\begin{equation} \label{1.14-3}
 |\partial_x^k( \hat{\rho}-\rho_+, \hat{u}-u_+, \hat{\theta}-\theta_+)|\leq
 C\delta\exp(-\hat{c}x),\quad k=0,1,2,
\end{equation}where $\delta=|(u_--u_+,\theta_--\theta_+)|$ and $C,
\hat{c}$ are positive constants independent of $x, \delta$.

2) For subsonic case $M_+<1$, there exists a positive constant
$\delta_0$ and a certain subset $\mathcal {M}^-\subset \mathcal
{M}_{\delta_0}$ such that if $(u_-, \theta_-)\in \mathcal {M}^-$,
then there exists a unique solution $(\hat{\rho}, \hat{u},
\hat{\theta})(x)$ satisfying \eqref{1.14-3}.

3) For transonic case  $M_+=1$, under the additional assumption
\begin{equation}\label{1.15-3}
\begin{aligned}
 &p_\theta(v_+,\theta_+)>0,\,\,p_{vv}(v_+,\theta_+)\geq 0
 ,\,\,p_{\theta\theta}(v_+,\theta_+)\geq 0,\\
 &p_{v\theta}(v_+,\theta_+)\leq 0
 ,\,\,e_{vv}(v_+,\theta_+)\leq 0,\,\,
 e_{\theta\theta}(v_+,\theta_+)\leq 0,
\end{aligned}
\end{equation}
there exists a positive constant $\delta_0$ and a certain curve
$\mathcal {M}^0\subset \mathcal {M}_{\delta_0}$ such that if $(u_-,
\theta_-)\in \mathcal {M}^0$, then there exists a unique solution
$(\hat{\rho}, \hat{u}, \hat{\theta})(x)$ satisfying
\begin{equation} \label{1.14-3-3}
 |\partial_x^k( \hat{\rho}-\rho_+, \hat{u}-u_+, \hat{\theta}-\theta_+)|\leq
 \frac{C\delta^{k+1}}{(1+\delta x)^{k+1}}+C\delta\exp(-\hat{c}x),\quad k=0,1,2,
\end{equation}
 and
\begin{equation} \label{1.15-3-3}(\hat{u}_{x}, \hat{\theta}_{x})
= (a_1,a_2)z^2(x)+O(z^3(x)+\delta\exp(-cx)).
\end{equation}where $a_i(i=1,2)$ are positive constants depending
only on $\mu,\kappa,\rho_+,u_+,\theta_+$ and $z(x)$ is a smooth
function satisfying
\begin{equation} \label{1.15-33}
0<c\frac{\delta}{1+\delta x}\leq z(x)\leq C\frac{\delta}{1+\delta
x},\quad |\partial_x^kz(x)|\leq C\frac{\delta^{k+1}}{(1+\delta
x)^{k+1}},\,\,k=1,2.
\end{equation}
\end{proposition}

\begin{rmk}
Note that the condition \eqref{1.15-3} holds when the gas is ideal
polytropic at $x=+\infty$.
\end{rmk}

Our main result in this paper is the following Theorem:
\begin{theorem} \label{theo 1.1}(Stability of transonic stationary solution).
Let $\rho_+>0,u_\pm<0,\theta_\pm>0$.  Assume that  $(\hat{\rho},
\hat{u}, \hat{\theta})(x)$ is the solution of the system
\eqref{1.12-3} and \eqref{1.12-31} satisfying \eqref{1.14-3-3} and
\eqref{1.15-3-3} for transonic case  $M_+=1$, and
\begin{equation}\label{1.17-3-3}
\begin{aligned}
&\beta_1:=p^+_{vv}+\frac{p^+_v}{v_+}>0,\quad
\beta_2:=\beta_1 p_{ss}^+-(p^+_{vs}+\frac{p^+_s}{2 v_+})^2>0,\\
& \beta_3:=-4v_+p^+_{v}\beta_2-(p_{ s}^+)^2\beta_1>0,\end{aligned}
\end{equation}where $p^+=p(v_+, s_+)$.
Also, suppose that the initial data $(\rho_0,u_0,\theta_0)$
satisfies
\begin{equation}\label{1.17-3}
(\rho_0-\hat{\rho},u_0-\hat{u},\theta_0-\hat{\theta})\in H^{1}(0,
\infty), \quad u_0(0)=u_-,\,\,\, \theta_0(0)=\theta_-.
\end{equation} Then, there exists a positive
constant $\varepsilon_{0}$ such that if
\begin{equation}\label{1.18-3}
\|(\rho_0-\hat{\rho},u_0-\hat{u},\theta_0-\hat{\theta})\|_1+\delta\leq
\varepsilon_0,\end{equation} where $\delta=|(u_--u_+,
\theta_--\theta_+)|$,  the outflow problem
\eqref{1.1}-\eqref{1.3} has a unique global solution
$(\rho,u,\theta)(x,t)$ satisfying $$\begin{aligned}
&(\rho-\hat{\rho},u-\hat{u},\theta-\hat{\theta})\in C([0, \infty); H^{1}(0, \infty)),\\
&\rho_{x} \in L _{2}(0, \infty;L_{2}(0, \infty)),\,\,\,
u_{x},\theta_{x}  \in L _{2}(0, \infty;H^{1}(0, \infty)).
\end{aligned}$$
Moreover, the solution $(\rho,u,\theta)(x,t)$ convergence to the
stationary solution $(\hat{\rho}, \hat{u}, \hat{\theta})(x)$
uniformly as time tends to infinity:
 \begin{equation} \label{1.19-3}
\lim_{t\rightarrow\infty}\sup_{x\in(0,\,\infty)}|(\rho, u,
\theta)(x,t)-(\hat{\rho}, \hat{u}, \hat{\theta})(x)|
 =0.
\end{equation}
\end{theorem}

\begin{rmk}The condition \eqref{1.17-3-3} plays the essential
role in the energy estimate for the transonic case (see Lemma
\ref{lem 3.2}). Note that the condition \eqref{1.17-3-3} holds when
the gas is ideal polytropic at $x=+\infty$. In fact, using
\eqref{perfect}, it is easy
 to check that
$$\beta_1=\gamma^2 v_+^{-2}p^+,\quad
\beta_2=\f{(\gamma-1)^2}{R^2}(\gamma-\f{1}{4}) v_+^{-2}(p^+)^2,\quad
\beta_3=\f{(\gamma-1)^2}{R^2}\gamma(3\gamma-1) v_+^{-2}(p^+)^3,$$ for
ideal polytropic gas.
\end{rmk}

The outline of this paper is organized as follows. Section 2 is
devoted to Theorem \ref{theo 1.1}. In Appendix, we will give the proof of
Proposition \ref{theo 1.1-3}.

\section{Stability of transonic stationary solution}
\setcounter{equation}{0}
In this section, we will give the proof of Theorem 1.1. We rewrite
\eqref{1.1} and \eqref{1.12-3} respectively as
\begin{equation}\label{3.1-3}\begin{cases}
\rho_{t}+(\rho u)_{x}= 0, \qquad x>0,\,\,\, t>0,\\
\rho (u_{t}+uu_x)+p_{x}
=\mu u_{xx},\\
\rho (e_{t}+ue_x)+pu_{x}
=\kappa\theta_{xx}+\mu u_{x}^2,\\
\rho \theta(s_{t}+us_x)=\kappa\theta_{xx}+\mu u_{x}^2
\end{cases}
\end{equation} and \begin{equation}\label{3.2-3}\begin{cases}
(\hat{\rho} \hat{u})_{x}= 0, \qquad x>0,\,\,\, t>0,\\
\hat{\rho} \hat{u}\hat{u}_x+\hat{p}_{x}
=\mu \hat{u}_{xx},\\
\hat{\rho} \hat{u}\hat{e}_x+\hat{p}\hat{u}_{x}
=\kappa\hat{\theta}_{xx}+\mu \hat{u}_{x}^2,\\
\hat{\rho}\hat{\theta} \hat{u}\hat{s}_x=\kappa\hat{\theta}_{xx}+\mu
\hat{u}_{x}^2,
\end{cases}
\end{equation}where $s=s(\rho,\theta)$ denotes the entropy and $\hat{p}=p(\hat{\rho},\hat{\theta}),\,\,
\hat{e}=e(\hat{\rho},\hat{\theta}),\,\,\hat{s}=s(\hat{\rho},\hat{\theta})$.
We set the perturbation $(\varphi,\psi,\zeta)(x,t)$ by
$$
 (\varphi,\psi,\zeta)(x,t)=(\rho,u,\theta)(x,t)-
 (\hat{\rho},\hat{u},\hat{\theta})(x),
$$ and
the solution space $ X(\emph{I})$ as
$$
X(\emph{I})=\{(\varphi,\psi,\zeta)\in C(\emph{I}; H^{1})\mid
\varphi_{x} \in L _{2}(\emph{I};L_{2}), (\psi_{x},\zeta_{x}) \in L
_{2}(\emph{I};H^{1})\},$$ for any interval $\emph{I}\subset
[0,\infty)$.
To prove Theorem \ref{theo 1.1}, it is sufficient to show the
following a priori estimate
\begin{proposition}\label{proposition 3.2} (A priori
 estimate) Besides the assumptions of Theorem \ref{theo 1.1}, suppose  that
$(\rho,u,\theta)$ is the solution to the outflow problem
\eqref{1.1}-\eqref{1.3} satisfying $(\phi,\psi,\zeta)\in X ([0,T])$.
Then, there exists a positive constant $\varepsilon_{1}$  such that
if $\,\sup_{0\leq t \leq T}\|(\varphi, \psi,\zeta)(t)\|_1\leq
\varepsilon_{1}$ and  $\delta=|(u_--u_+,\theta_--\theta_+)|\leq
\varepsilon_{1}$, then for any $t \in [0, T]$, it holds that
\begin{equation}\label{3.3-3}\begin{aligned}
 \|(\varphi, \psi, \zeta)(t)\|_1^{2}&+\int_{0}^{t}\lf(\|\varphi_{x}\|^{2}+
 \|(\psi_{x},\zeta_{x})(\tau)\|_1^{2}\rg)d\tau
 \\
 &+\int_{0}^{t}\lf(|\varphi(0,\tau)|^{2}+
 |\varphi_{x}(0,\tau)|^{2}\rg)d\tau
 \leq C\|(\varphi, \psi, \zeta)(0)\|_1^{2},
 \end{aligned}\end{equation}where $C$ is a positive constant independently of
 $t, T, \varepsilon_1$.
\end{proposition}

To prove Proposition \ref{proposition 3.2}, a crucial step is the
following energy estimate:
\begin{lemma}\label{lem 3.2} Under the assumptions of Proposition \ref{proposition 3.2}, it follows that
\begin{equation}\label{3.4-3}\begin{aligned}
\|(\varphi,\psi,\zeta)(t)\|^2&+\int_{0}^{t}\|(\psi_x,\zeta_x)(\tau)\|^2
d\tau+ \int_{0}^{t}|\varphi(0,\tau)|^2 d\tau\\
&\leq
C\|(\varphi,\psi,\zeta)(0)\|^2+C\delta\int_{0}^{t}\|\varphi_x(\tau)\|^2
d\tau.
\end{aligned}\end{equation}
\end{lemma}
\textbf{Proof}.
 Let
\begin{equation}\label{enery form}\begin{aligned}
\mathcal{E}:&=(e-\hat{e})-\hat{\theta}(s-\hat{s})+\f{\psi^2}{2}+\hat{p}\lf(\frac{1}{\rho}-\frac{1}{\hat{\rho}}\rg)\\
&=(e-\hat{\theta}s)+\f{\psi^2}{2}+\hat{p}\lf(v-\hat{v}\rg)-(\hat{e}-\hat{\theta}\hat{s}),
\end{aligned}\end{equation} then by \eqref{3.1-3} and \eqref{3.2-3}, we have
\begin{equation}\label{3.5-3}\begin{aligned}
(\hat{\rho}\mathcal{E})_t&+(\hat{\rho}
\hat{u}\mathcal{E})_x=\hat{\rho}\mathcal{E}_t+(\hat{\rho}
\hat{u})\mathcal{E}_x\\
&=\frac{\hat{\rho}}{\rho}\lf(1-\frac{\hat{\theta}}{\theta}\rg)(\kappa\theta_{xx}+\mu
u_{x}^2) -\frac{\hat{\rho}}{\rho}pu_x-\hat{\rho}
\hat{u}s\hat{\theta}_x+(\hat{u}-u)(\hat{\rho}e_x-\hat{\rho}\hat{\theta}s_x)\\
&+\frac{\hat{\rho}}{\rho}\mu\psi\psi_{xx}+\mu\hat{\rho}\psi\hat{u}_{xx}(v-\hat{v})-\hat{\rho}\psi
\lf(\frac{p_x}{\rho}-\frac{\hat{p}_x}{\hat{\rho}}\rg)-\hat{\rho}\psi^2\hat{u}_{x}-\hat{\rho}\psi^2\psi_{x}\\
&+\hat{\rho}
\hat{u}\lf(\hat{p}_x(v-\hat{v})-\hat{p}\hat{v}_x\rg)+\frac{\hat{\rho}}{\rho}\hat{p}u_x
+\hat{\rho}\hat{p}(\hat{u}-u)v_x\\
&-\hat{\rho} \hat{u}\lf(\hat{e}_x-
\hat{\theta}_x\hat{s}-\hat{\theta}\hat{s}_x\rg),
\end{aligned}\end{equation}where $\hat{v}=\hat{\rho}^{-1}.$

Arranging \eqref{3.5-3} yields
\begin{equation}\label{3.5-333}
\begin{aligned} (\hat{\rho}\mathcal{E})_t&+(\hat{\rho}
\hat{u}\mathcal{E})_x=\hat{\rho}\mathcal{E}_t+(\hat{\rho}
\hat{u})\mathcal{E}_x\\
&=\frac{\hat{\rho}}{\rho}\lf(1-\frac{\hat{\theta}}{\theta}\rg)(\kappa\theta_{xx}+\mu u_{x}^2)\\
&-\hat{\rho}\hat{ u}s\hat{\theta}_x+\hat{\rho}\hat{
u}\lf(\hat{p}_x(v-\hat{v})-\hat{p}\hat{v}_x\rg)-\hat{\rho}\hat{
u}\lf(\hat{e}_x-
\hat{\theta}_x\hat{s}-\hat{\theta}\hat{s}_x\rg)\\
&-\frac{\hat{\rho}}{\rho}pu_x+\frac{\hat{\rho}}{\rho}\hat{p}u_x-\hat{\rho}\psi
\lf(\frac{p_x}{\rho}-\frac{\hat{p}_x}{\hat{\rho}}\rg)\\
&+\mu\frac{\hat{\rho}}{\rho}\psi\psi_{xx}+\mu\hat{\rho}\psi\hat{u}_{xx}(v-\hat{v})-\hat{\rho}\psi^2\hat{u}_{x}
-\hat{\rho}\psi^2\psi_{x}\\
&+(\hat{u}-u)(\hat{\rho}e_x-\hat{\rho}\hat{\theta}s_x)+\hat{\rho}\hat{p}(\hat{u}-u)v_x\\
&=I_1+I_2+I_3+I_4+I_5.
\end{aligned}\end{equation}
Noticing that
$$\begin{aligned}\frac{\rho}{\hat{\rho}}I_1=&\kappa\lf(\f{\zeta\zeta_x}{\theta
}\rg)_x-\kappa\f{\hat{\theta}}{\theta^2}\zeta_x^2+\kappa\f{\hat{\theta}_x\zeta\zeta_x}{\theta^2}
-(\kappa\hat{\theta}_{xx}+\mu \hat{u}_{x}^2)\f{\zeta^2}{\theta
\hat{\theta}}\\
&+\mu\f{\zeta}{\theta}\lf(\psi_x^2+2\psi_x\hat{u}_x
\rg)+\hat{\rho}\hat{u} \hat{s}_{x}(\theta-
\hat{\theta}),\end{aligned}$$
we have
$$\begin{aligned}I_1=&\kappa\lf(\f{\hat{\rho}\zeta\zeta_x}{\rho\theta
}\rg)_x-\kappa\lf(\f{\varphi}{\rho}\rg)_x\f{\zeta\zeta_x}{\theta
}-\kappa\f{\hat{\rho}\hat{\theta}}{\rho\theta^2}\zeta_x^2+\kappa\f{\hat{\rho}\hat{\theta}_x\zeta\zeta_x}{\rho\theta^2}
\\
&-(\kappa\hat{\theta}_{xx}+\mu
\hat{u}_{x}^2)\f{\hat{\rho}\zeta^2}{\rho\theta
\hat{\theta}}+\mu\f{\hat{\rho}\zeta}{\rho\theta}\lf(\psi_x^2+2\psi_x\hat{u}_x\rg)
+\frac{\hat{\rho}}{\rho}\hat{\rho}\hat{u} \hat{s}_{x}(\theta-
\hat{\theta}).\end{aligned}$$
Using
$e_x=-pv_x+\theta s_x$, we have
$$\begin{aligned}&I_2=\hat{\rho}
\hat{u}\lf(\hat{p}_x(v-\hat{v})-\hat{\theta}_x(s-\hat{s})\rg),\\
&I_5=\hat{\rho}
(\hat{u}-u)\lf((\theta-\hat{\theta})s_x-(p-\hat{p})v_x\rg).\end{aligned}$$
Also, we can arrange as follows:
$$I_3
=-\lf(\frac{\hat{\rho}}{\rho}(p-\hat{p})\psi\rg)_x+\lf(\frac{\varphi}{\rho}\rg)_x(p-\hat{p})\psi
-\frac{\hat{\rho}}{\rho}(p-\hat{p})\hat{u}_x-
\hat{\rho}\psi\hat{p}_x(v-\hat{v})$$ and
$$I_4
=\mu\lf(\frac{\hat{\rho}}{\rho}\psi\psi_{x}\rg)_x-\mu\lf(\frac{\varphi}{\rho}\rg)_x\psi\psi_{x}
-\mu\frac{\hat{\rho}}{\rho}\psi_{x}^2+\mu\hat{\rho}\psi\hat{u}_{xx}(v-\hat{v})-\hat{\rho}\psi^2\hat{u}_{x}
-\hat{\rho}\psi^2\psi_{x}.$$
Substituting $I_i(i=1,\cdots,5)$ into \eqref{3.5-3}, we have
\begin{equation}\label{3.9-3}
(\hat{\rho}\mathcal{E})_t+(\hat{\rho}
\hat{u}\mathcal{E})_x+\mu\f{\hat{\rho}\hat{\theta}}{\rho\theta}\psi_x^2
+\kappa\f{\hat{\rho}\hat{\theta}}{\rho\theta^2}\zeta_x^2=
\Delta_{1x}+\Delta_2+\Delta_3+\Delta_4,
\end{equation}where
$$\begin{aligned}
&\Delta_{1}=\frac{\mu\hat{\rho}\psi\psi_x}{\rho}+\kappa\f{\hat{\rho}\zeta\zeta_x}{\rho\theta
}-\frac{\hat{\rho}(p-\hat{p})\psi}{\rho}-\frac{\hat{\rho}\psi^3}{3},\\
&\Delta_{2}=
\kappa\f{\hat{\rho}\hat{\theta}_x\zeta\zeta_x}{\rho\theta^2}
+2\mu\f{\hat{\rho}\zeta}{\rho\theta}\psi_x\hat{u}_x
+\mu\hat{\rho}\psi\hat{u}_{xx}(v-\hat{v})
+\frac{\hat{\rho}_x\psi^3}{3},\\
&\Delta_{3}=-\kappa\lf(\f{\varphi}{\rho}\rg)_x\f{\zeta\zeta_x}{\theta
}+\lf(\frac{\varphi}{\rho}\rg)_x(p-\hat{p})\psi-\mu\lf(\frac{\varphi}{\rho}\rg)_x\psi\psi_{x},\\
&\begin{aligned}\Delta_{4}=&-(\kappa\hat{\theta}_{xx}+\mu
\hat{u}_{x}^2)\f{\hat{\rho}\zeta^2}{\rho\theta
\hat{\theta}}-\hat{\rho}\psi^2\hat{u}_{x}-\frac{\hat{\rho}}{\rho}(p-\hat{p})\hat{u}_x-
\hat{\rho}\psi\hat{p}_x(v-\hat{v})\\
&+\hat{\rho}
\hat{u}\lf(\hat{p}_x(v-\hat{v})-\hat{\theta}_x(s-\hat{s})\rg)
+\frac{\hat{\rho}}{\rho}\hat{\rho}\hat{u}\hat{s}_{x}(\theta-\hat{\theta})\\
&+\hat{\rho}
(\hat{u}-u)\lf((\theta-\hat{\theta})s_x-(p-\hat{p})v_x\rg).
\end{aligned}\end{aligned}$$
It is easy to check that there exist positive constants $c_i(i=1,2)$
satisfying
\begin{equation}\label{3.10-3}
c_1(\varphi^2+\psi^2+\zeta^2)\leq \mathcal{E}(x,t)\leq
c_2(\varphi^2+\psi^2+\zeta^2),\end{equation} due to \eqref{1.10-3}
and the assumptions of Proposition \ref{proposition 3.2}.
By using \eqref{3.10-3} and  $u\mid_{x=0}=u_-<0,\,\,\,(\psi,
\zeta)\mid_{x=0}=0$,  we have
\begin{equation}\label{3.11-3}
\Delta_1\mid_{x=0}=0,\quad - (\rho u\mathcal{E})\mid_{x=0}\geq
c_3\varphi^2(0, t).\end{equation}
Integrating \eqref{3.9-3} for  $(x,t)$ and using \eqref{3.10-3} and
\eqref{3.11-3}, we have
\begin{equation}\label{3.12-3}\begin{aligned}
\|(\varphi,\psi,\zeta)(t)\|^2&+\int_{0}^{t}\|(\psi_x,\zeta_x)(\tau)\|^2
d\tau+ \int_{0}^{t}|\varphi(0,\tau)|^2 d\tau\\
&\leq
C\|(\varphi,\psi,\zeta)(0)\|^2+C\sum_{i=2}^4\int_{0}^{t}\int_{0}^{\infty}|\Delta_{i}|dxd\tau.
\end{aligned}\end{equation}
By \eqref{1.14-3-3}, we have
$$|\Delta_{2}|\leq
C\delta|(\psi_x,\zeta_{x})|^2+C\delta\lf(\frac{\delta^2}{(1+\delta
x)^4}+e^{-\hat{c}x}\rg)|(\varphi,\psi,\zeta)|^2
+C\frac{\delta^2}{(1+\delta x)^2}|\psi|^3$$ and using the inequality
\begin{equation} \label{ineq-3}
|f(x)|\leq |f(0)|+\sqrt{x}\|f_x\|,\,\,\forall f\in H^1(0, \infty)
\end{equation}yields
\begin{equation}\label{3.13-3}\int_{0}^{\infty}|\Delta_{2}|dx\leq
C\delta\|(\varphi_x, \psi_x,\zeta_x))\|^2+
C\delta|\varphi(0,t)|^2+C\varepsilon_1\delta^\frac{1}{2}\|\psi_x\|^2,\end{equation}
where we used $$\begin{aligned}
\int_{0}^{\infty}&\frac{\delta^2}{(1+\delta x)^2}|\psi|^3dx\leq
C\|\psi\|\|\psi_x\|
\int_{0}^{\infty}\frac{\delta^2}{(1+\delta x)^2}|\psi|dx\\
&\leq C\|\psi\|\|\psi_x\|^2
\int_{0}^{\infty}\frac{\delta^2\sqrt{x}}{(1+\delta x)^2}dx\leq
C\varepsilon_1\delta^\frac{1}{2}\|\psi_x\|^2.
\end{aligned}$$
Noticing that
$$|\Delta_{3}|\leq
C(|\varphi_x|+|\hat{\rho}_x||\varphi|)(|\zeta||\zeta_x|+|(\varphi,\psi,\zeta)|^2+|\psi||\psi_x|)$$
and by the  argument similar to \eqref{3.13-3}, we have
$$\int_{0}^{\infty}|\Delta_{3}|dx\leq
C(\delta+\varepsilon_1)\|(\varphi_x, \psi_x,\zeta_x)\|^2+
C(\delta+\varepsilon_1)|\varphi(0,t)|^2.$$
It is more difficult to
estimate $\Delta_{4}$. We first rewrite $\Delta_{4}$ as
\begin{equation}\label{3.14.33}\begin{aligned}\Delta_{4}=&-(p-\hat{p})\hat{u}_x+
\hat{\rho}\hat{u}\hat{p}_x(v-\hat{v})-\hat{\rho}
\hat{u}\hat{\theta}_x(s-\hat{s})
+\hat{\rho}\hat{u}\hat{s}_{x}(\theta-\hat{\theta})\\&-(\kappa\hat{\theta}_{xx}+\mu
\hat{u}_{x}^2)\f{\hat{\rho}\zeta^2}{\rho\theta
\hat{\theta}}-\hat{\rho}\psi^2\hat{u}_{x}\\
&+\hat{\rho}
(\hat{u}-u)\lf((\theta-\hat{\theta})s_x-(p-\hat{p})v_x+\hat{p}_x(v-\hat{v})\rg)\\
&+\lf(1-\frac{\hat{\rho}}{\rho}\rg)\lf(\hat{u}_x(p-\hat{p})-
\hat{\rho}\hat{u}\hat{s}_{x}(\theta-\hat{\theta})\rg)=:\sum_{i=1}^4\Delta_{4}^i.
\end{aligned}\end{equation}
Using $$-\theta_v(v,s)=p_s(v,s),\,\,\,
\hat{p}_x=p_v(\hat{v},\hat{s})\hat{v}_x+p_s(\hat{v},\hat{s})\hat{s}_x,
\,\,\, \hat{u}\hat{v}_x=\hat{u}_x\hat{v}$$ and
$\hat{\rho}\hat{\theta} \hat{u}\hat{s}_x=\kappa\hat{\theta}_{xx}+\mu
\hat{u}_{x}^2$ yields
$$\begin{aligned}&\begin{aligned}
\Delta_{4}^1=&-\hat{u}_x\lf(p-\hat{p}-\hat{p}_v(v-\hat{v})-\hat{p}_s(s-\hat{s})\rg)\\
&+\hat{\rho}\hat{u}\hat{s}_x\lf(\theta-\hat{\theta}-\hat{\theta}_v(v-\hat{v})-\hat{\theta}_s(s-\hat{s})\rg)
,\end{aligned}\\&\Delta_{4}^2=-\hat{\rho}\psi^2\hat{u}_{x}-\frac{\hat{\rho}^2\hat{u}\hat{s}_x}{\rho\theta}\zeta^2,
\\&\begin{aligned}
\Delta_{4}^3=&\hat{\rho}
(\hat{u}-u)\lf(\hat{p}_x(v-\hat{v})-(p-\hat{p})\hat{v}_x+(\theta-\hat{\theta})\hat{s}_x\rg)\\&+\hat{\rho}
(\hat{u}-u)\lf((\theta-\hat{\theta})(s-\hat{s})_x-(p-\hat{p})(v-\hat{v})_x\rg),
\end{aligned}\end{aligned}$$
where
$\hat{p}_v=p_v(\hat{v},\hat{s}),\,\,\,\hat{p}_s=p_s(\hat{v},\hat{s}),
\,\,\,\hat{\theta}_v=\theta_v(\hat{v},\hat{s}),\,\,\,\hat{\theta}_s=\theta_s(\hat{v},\hat{s}).$
Substituting $\Delta_{4}^i(i=1,\cdots,4)$ into \eqref{3.14.33} and
using
$$\begin{aligned}
&
\hat{p}_x=p_v(\hat{v},\hat{s})\hat{v}_x+p_s(\hat{v},\hat{s})\hat{s}_x,\\
&p-\hat{p}=p_v(\hat{v},\hat{s})(v-\hat{v})+p_s(\hat{v},\hat{s})(s-\hat{s})
+O((v-\hat{v})^2+(s-\hat{s})^2),\end{aligned}$$ we have
\begin{equation}\label{3.14.333}
\begin{aligned}
\Delta_{4}=&-\hat{u}_x\lf(p-\hat{p}-\hat{p}_v\phi-\hat{p}_s\chi\rg)-\hat{\rho}\hat{u}_x\psi^2
-\hat{\rho}\hat{p}_s\hat{v}_x\psi\chi-\hat{\rho}\hat{p}_v
\hat{u}_x\phi^2-\hat{\rho}\hat{p}_s \hat{u}_x\phi\chi
\\
&+\hat{\rho}\hat{u}\hat{s}_x\lf(\theta-\hat{\theta}-\hat{\theta}_v(v-\hat{v})-\hat{\theta}_s(s-\hat{s})\rg)
-\hat{\rho}^2\hat{u}\hat{s}_x\lf(\frac{\zeta^2}{\rho
u}+\phi\zeta\rg)\\
&+\hat{\rho}\psi \lf(\zeta\hat{s}_x+\hat{p}_s\hat{s}_x\phi+\zeta\chi_x-(p-\hat{p})\phi_x\rg)
+O(\phi^2+\chi^2)\hat{\rho}(\psi\hat{v}_x+\phi\hat{u}_x)\\
=&J_1+J_2+J_3,
\end{aligned}
\end{equation}
where $\phi=v-\hat{v},\,\,  \chi=s-\hat{s}$.
By the argument similar to \eqref{3.13-3}, we have
\begin{equation}\label{3.15-33}
\int_{0}^{\infty}(|J_{2}|+|J_{3}|)dx\leq
C(\delta+\varepsilon_1)\|(\varphi_x, \psi_x,\zeta_x))\|^2+
C(\delta+\varepsilon_1)|\varphi(0,t)|^2,
\end{equation}
where we used the fact that $|\hat{s}_x|\leq
C(|\hat{\theta}_{xx}|+|\hat{u}_x|^2)\leq C\frac{\delta^3}{(1+\delta
x)^3}$.
By using $\hat{\rho}\hat{u}\hat{v}_x=\hat{u}_x$ and
$$
p-\hat{p}-\hat{p}_v\phi-\hat{p}_s\chi=\hat{p}_{vv}\phi^2+2\hat{p}_{vs}\phi\chi+\hat{p}_{ss}\chi^2+O(\phi^3+\chi^3),
$$
we have
$$\begin{aligned}
J_1=&-\hat{u}_x\lf((p^+_{vv}+\frac{p^+_v}{v_+})\phi^2+p^+_{ss}\chi^2
+\rho_+\psi^2 +(2p^+_{vs}
+\frac{p^+_s}{v_+})\phi\chi+\frac{p^+_s}{u_+}\psi\chi\rg)\\
&-\hat{u}_x\lf(\lf((\hat{p}_{vv}-p^+_{vv})+(\frac{\hat{p}_v}{\hat{v}}
-\frac{p^+_v}{v_+})\rg)\phi^2
+(\hat{p}_{ss}-p^+_{ss})\chi^2+(\hat{\rho}-\rho_+)\psi^2
\rg)\\
&-\hat{u}_x\lf(\lf(2(\hat{p}_{vs}-p^+_{vs})
+(\frac{\hat{p}_s}{\hat{v}}-\frac{p^+_s}{v_+})\rg)\phi\chi
+(\frac{\hat{p}_s}{\hat{u}}-\frac{p^+_s}{u_+})\psi\chi+O(\phi^3+\chi^3)\rg)\\
=&:-\sum_{i=1}^3J_1^i,
\end{aligned}
$$
where
$p^+=p(v_+,s_+), \,\hat{p}=p(\hat{v},\hat{s}).$
By the  argument similar to \eqref{3.13-3}, we have
\begin{equation}\label{3.14-3}
\int_{0}^{\infty}(|J_{1}^2|+|J_1^{3}|)dx\leq
C(\delta+\varepsilon_1)\|(\varphi_x, \psi_x,\zeta_x))\|^2+
C(\delta+\varepsilon_1)|\varphi(0,t)|^2.
\end{equation}
To estimate $J_1^1$, we define the quadratic form $f(\phi, \chi,
\psi)$ by
$$
f(\phi, \chi,
\psi):=(p^+_{vv}+\frac{p^+_v}{v_+})\phi^2+p^+_{ss}\chi^2
+\rho_+\psi^2 +(2p^+_{vs}
+\frac{p^+_s}{v_+})\phi\chi+\frac{p^+_s}{u_+}\psi\chi.
$$
The matrix $A$ corresponding to the quadratic form $f$
$$
A:=\lf(\begin{matrix} (p^+_{vv}+\frac{p^+_v}{v_+}) &
(p^+_{vs}+\frac{p^+_s}{2 v_+})
& 0\\
(p^+_{vs}+\frac{p^+_s}{2 v_+}) & p^+_{ss} & \frac{p^+_s}{2 u_+}\\
 0& \frac{p^+_s}{2 u_+}  & \rho_+
\end{matrix}\rg),$$
is positive if and only if all principal minors $\bar{A}_i(i=1,2,3)$
of $A$ are positive.
Noticing that $$M_+=1 \rightleftarrows -p^+_{v}=(\rho_+ u_+)^2,
$$ and by using  \eqref{1.17-3-3}, we compute
$\bar{A}_i(i=1,2,3)$ as follows:
\begin{equation}\label{3.17-3}\begin{aligned}
&\bar{A}_1=p^+_{vv}+\frac{p^+_v}{v_+}>0,\\
&
 \bar{A}_2=(p^+_{vv}+\frac{p^+_v}{v_+}) p_{ss}^+-(p^+_{vs}+\frac{p^+_s}{2 v_+})^2>0,\\
& \begin{aligned}\bar{A}_3&=\det
A=\rho_+\bar{A}_2-\frac{(p_{s}^+)^2}{(2
u_+)^2}\bar{A}_1\\
&=\frac{1}{(2 u_+)^2}\lf(-4v_+p^+_{v}\bar{A}_2-(p_{
s}^+)^2\bar{A}_1\rg)>0.\end{aligned}
\end{aligned}\end{equation}
Using \eqref{3.17-3}, \eqref{1.15-3-3} and \eqref{1.15-33}, we have
\begin{equation}\label{3.16-33}\begin{aligned}
J_1^1=&-\hat{u}_xf(\phi, \chi, \psi)\\
&\leq- c_4z(x)^2|(\phi,\psi,\chi)|^2+C\lf(\f{\delta^3}{(1+\delta
x)^3}+\delta e^{-\hat{c}x}\rg)|(\phi,\psi,\chi)|^2
.\end{aligned}\end{equation}
Using \eqref{3.16-33}, we have by the same lines as in
\eqref{3.13-3}
\begin{equation}\label{3.16-3}
\int_{0}^{\infty}(|J_{1}^1|+z(x)^2|(\phi,\psi,\chi)|^2)dx\leq
C\delta\|(\varphi_x, \psi_x,\zeta_x))\|^2+ C\delta|\varphi(0,t)|^2.
\end{equation}
By \eqref{3.15-33}, \eqref{3.14-3} and \eqref{3.16-3}, we have
$$\int_{0}^{\infty}|\Delta_{4}|dx\leq
C(\delta+\varepsilon_1)\|(\varphi_x, \psi_x,\zeta_x))\|^2+
C(\delta+\varepsilon_1)|\varphi(0,t)|^2.$$
Substituting the estimates for $\Delta_i(i=2,3,4)$ into
\eqref{3.12-3}, and choosing $\delta$ to be small, we obtain
\eqref{3.4-3}.The proof of Lemma \ref{lem 3.2} is completed.


\begin{lemma}\label{lem 3.3} Under the assumptions of Proposition \ref{proposition 3.2}, it follows that
\begin{equation}\label{3.18-3}
\|(\varphi_x, \psi_x, \zeta_x)(t)\|^2+\int_{0}^{t}\|(\varphi_x,
\psi_{xx}, \zeta_{xx})(\tau)\|^2 d\tau+
\int_{0}^{t}|\varphi_x(0,\tau)|^2 d\tau\leq
C\|(\varphi,\psi,\zeta)(0)\|_1^2.
\end{equation}
\end{lemma}
\textbf{Proof}. We first estimate $\|\varphi_x(t)\|$.
 Subtracting the first equation in \eqref{3.2-3} from the first equation in \eqref{3.1-3}
 and applying $\partial_x$ to the resulting equality yields
$$
\varphi_{xt}+u\varphi_{xx}+\rho\psi_{xx}=-2u_x\varphi_x-2\hat{\rho}_x\psi_x-
\hat{u}_{xx}\varphi-\hat{\rho}_{xx}\psi.
$$
 Also, subtracting the second equation in \eqref{3.2-3} from the
second equation in \eqref{3.1-3}
 yields
$$
\rho\psi_{t}+\rho
u\psi_{x}+(p-\hat{p})_{x}=\mu\psi_{xx}+(\hat{u}\hat{\rho}- \rho
u)\hat{u}_x
$$
We multiply above two equations by $\mu\frac{\varphi_x}{\rho^3}$ and
$\frac{\varphi_x}{\rho^2}$, respectively, to discover (see \cite{CHS
17})
\begin{equation}\label{3.21-3}
\lf(\frac{\mu\varphi^2_x}{2\rho^3}+\frac{\varphi_x\psi}{\rho}\rg)_t+\lf(\frac{\mu
u\varphi^2_x}{2\rho^3}-\frac{\varphi_t}{\rho}\psi\rg)_x+
\frac{p_\rho(\rho,\theta)}{\rho^2}\varphi^2_x=F_1,
\end{equation}where
$$\begin{aligned}F_1:=&-\mu\frac{\varphi_x}{\rho^3}(\hat{\rho}_x\psi_x+
\hat{u}_{xx}\varphi+\hat{\rho}_{xx}\psi)\\
&+\frac{\varphi_x}{\rho^2}(\hat{u}\hat{\rho}- \rho
u)\hat{u}_x-\frac{u}{\rho}\varphi_x\psi_x+\frac{(\rho
u)_x\psi}{\rho^2}\hat{\rho}_x-\frac{(\rho u)_x}{\rho}\psi_x\\
&-\frac{\varphi_x}{\rho^2}\lf(p_\theta(\rho,\theta)\zeta_x
+\hat{\rho}_x\lf(p_\rho(\rho,\theta)-p_\rho(\hat{\rho},\hat{\theta})\rg)+
\hat{\theta}_x\lf(p_\theta(\rho,\theta)-p_\theta(\hat{\rho},\hat{\theta})\rg)\rg).\end{aligned}$$
By using \eqref{1.14-3-3},  $\psi\mid_{x=0}=0,\,\,u\mid_{x=0}=u_-<0$
and the assumptions of Proposition \ref{proposition 3.2}, we have
\begin{equation}\label{3.22-3}\begin{aligned}
&\int_0^\infty\lf(\frac{\mu
u\varphi^2_x}{2\rho^3}-\frac{\varphi_t}{\rho}\psi\rg)_xdx=-\f{\mu
u_-}{2\rho^3(0,t)}
\varphi_x^2(0,t)\geq c_5\varphi_x^2(0,t),\\
&\int_{0}^{\infty}|F_1|dx\leq C(\delta+\eta)\|\varphi_x\|^2+
C(1+\eta^{-1})\|(\psi_x,\zeta_x)\|^2+C\delta|\varphi(0,t)|^2,
\end{aligned}\end{equation} for any $\eta>0$.
After integrating \eqref{3.21-3} for  $(x,t)$,  and using
\eqref{3.22-3},  $p_\rho(\rho,\theta)>0$\ and \eqref{3.4-3}, we have
\begin{equation}\label{3.23-3}\|\varphi_x(t)\|^2+\int_{0}^{t}\|\varphi_x(\tau)\|^2 d\tau+
\int_{0}^{t}|\varphi_x(0,\tau)|^2 d\tau\leq
C\|(\varphi,\psi,\zeta,\varphi_x)(0)\|^2.\end{equation}
The following, we will estimate $\|(\psi_x, \zeta_x)\|$.
Subtracting the second equation in \eqref{3.2-3} from the second
equation in \eqref{3.1-3}
 and multiplying it by $-\frac{\psi_{xx}}{\rho}$ yields
\begin{equation}\label{3.25-3}\lf(\frac{\psi^2_x}{2}\rg)_t-(\psi_t\psi_x)_x+
\frac{\mu \psi^2_{xx}}{\rho}=F_2,\end{equation}where
$$F_2=u\psi_x\psi_{xx}+\frac{(p-\hat{p})_x}{\rho}\psi_{xx}+\frac{(\hat{\rho}
\hat{u}-\rho u)\hat{u}_x}{\rho}\psi_{xx}.$$
Also, subtracting the third equation in \eqref{3.2-3} from the third
equation in \eqref{3.1-3} and using
$e_t=e_\theta(\rho,\theta)\theta_t-e_\rho(\rho,\theta)(\rho u)_x$,
we have
$$\begin{aligned}
\rho e_\theta(\rho,\theta)\zeta_{t}&+\rho u\zeta_x-\kappa\zeta_{xx}
=\rho e_\rho(\rho,\theta)(\rho u-\hat{\rho} \hat{u})_x \\
&+(\hat{u}\hat{\rho}- \rho
u)\hat{e}_x-(pu_{x}-\hat{p}\hat{u}_{x})+\mu(u_x^2-\hat{u}_x^2).\end{aligned}
$$
Multiplying it by $-\frac{\zeta_{xx}}{\rho e_\theta(\rho,\theta)}$
yields
\begin{equation}\label{3.26-3}
\lf(\frac{\zeta^2_x}{2}\rg)_t-(\zeta_t\zeta_x)_x+
\frac{\kappa\zeta^2_{xx}}{\rho
e_\theta(\rho,\theta)}=F_3,\end{equation}where
$$\begin{aligned}F_3&=\frac{\zeta_{xx}}{\rho
e_\theta(\rho,\theta)}\lf(\rho u\zeta_x+(pu_{x}-\hat{p}\hat{u}_{x})
-\rho e_\rho(\rho,\theta)(\rho u-\hat{\rho} \hat{u})_x\rg) \\
&-\frac{\zeta_{xx}}{\rho
e_\theta(\rho,\theta)}\lf((\hat{u}\hat{\rho}- \rho
u)\hat{e}_x+\mu(u_x^2-\hat{u}_x^2)\rg).
\end{aligned}$$
Adding \eqref{3.25-3} and \eqref{3.26-3}, we get
\begin{equation}\label{3.27-3}
\frac{1}{2}\lf(\psi^2_x+\zeta^2_x\rg)_t-(\psi_t\psi_x+\zeta_t\zeta_x)_x+
\lf(\frac{\mu \psi^2_{xx}}{\rho}+\frac{\kappa\zeta^2_{xx}}{\rho
e_\theta(\rho,\theta)}\rg)=F_2+F_3.\end{equation}
By using \eqref{1.14-3-3},  $\psi\mid_{x=0}=0,\,\,u\mid_{x=0}=u_-<0$
and the assumptions of Proposition \ref{proposition 3.2}, we have
\begin{equation}\label{3.28-3}\begin{aligned}
&\int_0^\infty\lf(\psi_t\psi_x+\zeta_t\zeta_x\rg)_xdx=0,\\
&\int_{0}^{\infty}|F_2|dx\leq C\eta\|\psi_{xx}\|^2+
C(\delta+\eta^{-1})\|(\varphi_x, \psi_x,\zeta_x)\|^2+C\delta|\varphi(0,t)|^2,\\
&\int_{0}^{\infty}|F_3|dx\leq C\eta\|\zeta_{xx}\|^2+
C(\delta+\eta^{-1})\|(\varphi_x,
\psi_x,\zeta_x)\|^2+C\delta|\varphi(0,t)|^2
\end{aligned}\end{equation} for any $\eta>0$.
After integrating \eqref{3.27-3} for  $(x,t)$, and using
\eqref{3.28-3}, $e_\theta(\rho,\theta)>0$, \eqref{3.4-3} and
\eqref{3.23-3}, we have
\begin{equation}\label{3.29-3}\|(\psi_x,
\zeta_x)(t)\|^2+\int_{0}^{t}\|(\psi_{xx},\zeta_{xx})(\tau)\|^2
d\tau\leq C\|(\varphi,\psi,\zeta)(0)\|^2_1.\end{equation}
By  \eqref{3.23-3} and \eqref{3.29-3}, we obtain \eqref{3.18-3}. The
proof of Lemma \ref{lem 3.3} is completed.

On the basis of the above  Lemmas \ref{lem 3.2}-\ref{lem 3.3}, we get the Proposition
\ref{proposition 3.2} at once.  Therefore,  the Theorem 1.1 is obtained.

\section{Appendix}
\setcounter{equation}{0}

This section is devoted to Proposition \ref{theo 1.1-3}, i.e. we
prove the existence of the solution to the stationary problem
\eqref{1.12-3} and \eqref{1.12-31}.
Integrating  \eqref{1.12-3} over $[x, \infty)$ yields
\begin{equation}\label{2.1-3}\left\{\begin{aligned}
&\hat{\rho} \hat{u}= \rho_+ u_+,\,\,x>0, \\
&\hat{\rho} \hat{u}^2+\hat{p}
=\mu \hat{u}_{x}+\rho_+ u_+^2+p_+,\\
&\hat{\rho}\hat{u}(\hat{e}+\frac{\hat{u}^{2}}{2}) +\hat{p}\hat{u}
=\kappa\hat{\theta}_{x}+\mu\hat{u}
\hat{u}_{x}+\rho_+u_+(e_++\frac{u_+^2}{2})+p_+u_+,
\end{aligned}\right.
\end{equation}
where $p_+=p(v_+, \theta_+),\,\,\,e_+=e(v_+, \theta_+)$.
Integrating  the first equation in \eqref{1.12-3} over $[0,x)$ and
using \eqref{1.12-31} yields
\begin{equation}\label{2.2-3}
\hat{\rho} \hat{u}= \rho(0) u_-,\,\,x>0.
\end{equation}
By \eqref{2.1-3} and \eqref{2.2-3}, we have \eqref{2.3-3}. So, if
$u_+\geq 0$, there is no solution to the system \eqref{1.12-3} and
\eqref{1.12-31}.
Assume that $u_+< 0$, substituting
$\hat{u}=\frac{u_+}{v_+}\hat{v}(v_+=\rho_+^{-1},\,\,\hat{v}=\hat{\rho}^{-1})$
into $(3.1)_2$ and $(3.1)_3$, we have
\begin{equation}\label{2.4-3}\begin{aligned}
&\hat{v}_{x}=\frac{u_+}{\mu v_+}(\hat{v}-v_+ )+\frac{v_+}{\mu
u_+}(p(\hat{v}, \hat{\theta})-p_+ )
=:g_1(\hat{v}, \hat{\theta}),\\
&\hat{\theta}_{x}=\frac{u_+}{\kappa v_+}(e(\hat{v},
\hat{\theta})-e_+ )-\frac{u_+^3}{2\kappa v_+^3}(\hat{v}-v_+
)^2+\frac{u_+}{\kappa v_+}p_+ (\hat{v}-v_+ )=:g_2(\hat{v},
\hat{\theta}).
\end{aligned}\end{equation}
Setting
$$W=(\hat{v}, \hat{\theta}),\,G(W)=(g_1(\hat{v}, \hat{\theta}), g_2(\hat{v}, \hat{\theta}))
,W^-=(v_b, \theta_-),\,\,W^+=(v_+, \theta_+),$$ we consider the
system below
\begin{equation}\label{2.5-3}\begin{aligned}
&W_{x}=G(W),\quad x>0\\
&W(0)=W^-,\quad \lim_{x\rightarrow\infty}W(x)=W_+,
\end{aligned}\end{equation}where $v_b=\frac{v_+}{u_+}u_-$.
The Jacobian matrix of $G$ at  $W_+$ is
\begin{equation}\label{2.6-3}
J_+=\lf(\begin{matrix} \frac{v_+}{\mu
u_+}\lf((\frac{u_+}{v_+})^2+p_v(v_+, \theta_+)\rg)
& \frac{v_+}{\mu u_+}p_\theta(v_+, \theta_+)\\
\frac{u_+}{\kappa v_+}\lf(e_v(v_+, \theta_+)+p(v_+, \theta_+)\rg) &
\frac{u_+}{\kappa v_+}e_\theta(v_+, \theta_+)
\end{matrix}\rg)\end{equation} and
\begin{equation}\label{2.7-3}
\det J_+= \frac{1}{\mu \kappa
}\lf(\frac{u_+^2}{v_+^2}+\widetilde{p}_v(v_+, s_+)\rg)e_\theta(v_+,
\theta_+).
\end{equation}due to \eqref{1.6-3} and \eqref{1.7-3}, where $s_+=s(v_+, \theta_+)$.
Using \eqref{1.9-3}  and \eqref{2.7-3}, we have
\begin{equation}\label{2.8-3}\begin{aligned}
M^+>1(=1,\,<1)&\leftrightarrows
\lf(\frac{u_+^2}{v_+^2}+\widetilde{p}_v(v_+,
s_+)\rg)>0(=0,\,<0)\\
&\leftrightarrows \det J_+>0(=0,\,<0).\end{aligned}
\end{equation}
The two eigenvalues $\lambda_1,\,\lambda_2$  of $J_+$ must be
satisfied:
\begin{equation}\label{2.9-3}
\lambda^2-b\lambda+\det J_+=0
\end{equation}
where
\begin{equation}\label{2.10-3}
b=\frac{v_+}{\mu u_+}\lf(\frac{u_+^2}{v_+^2}+\widetilde{p}_v(v_+,
s_+)+\frac{\theta_+(p_\theta(v_+, \theta_+))^2}{e_\theta(v_+,
\theta_+)}\rg)+ \frac{u_+}{\kappa v_+}e_\theta(v_+, \theta_+).
\end{equation}
By using \eqref{2.8-3}-\eqref{2.10-3}, \eqref{1.8-3} and $u_+<0$, we
have
\begin{equation}\label{2.9}\begin{aligned}
&\text{if}\quad M^+>1, \quad\text{then} \quad
b<0\quad\text{and}\quad \lambda_1< \lambda_2<0,
\\&\text{if}\quad M^+<1,\quad\text{then} \quad
\lambda_1<0< \lambda_2,
\end{aligned}
\end{equation}
where we used
$$b^2-4\det J_+\geq \frac{4u_+^2}{\mu \kappa
v_+^2}\theta_+(p_\theta(v_+, \theta_+))^2>0. $$
Thus, it is easy to prove the case 1) and 2) in Proposition
\ref{theo 1.1-3} by the same lines as in \cite{KNNZ 10}.
We prove the case 3) in Proposition \ref{theo 1.1-3}.
Due to \eqref{2.6-3}-\eqref{2.8-3}, for  transonic case  $M_+=1$, we
have
\begin{equation}\label{2.12-3} J_+=\lf(\begin{matrix}
\frac{v_+}{\mu u_+}\frac{\theta_+(p_\theta^+)^2}{e_\theta^+}
& \frac{v_+}{\mu u_+}p_\theta^+\\
\frac{u_+}{\kappa v_+}\theta_+p_\theta^+ & \frac{u_+}{\kappa
v_+}e_\theta^+
\end{matrix}\rg)\equiv\lf(\begin{matrix} a_{11}
& a_{12}\\
a_{21} & a_{22}
\end{matrix}\rg),\,\,\,\det J_+=0.
\end{equation} and the eigenvalues of $J_+$ are $\lambda_1=0$ and $\lambda_2=a_{11}+a_{22}<0$
of which corresponding eigenvectors are given by
\begin{equation}\label{2.10.3}
r_1=\lf(\begin{matrix} 1\vspace{0.3cm}\\
-\frac{a_{11}}{a_{12}}
\end{matrix}\rg)\equiv\lf(\begin{matrix} 1\vspace{0.3cm}\\
-\frac{\theta_+p_\theta^+}{e_\theta^+}
\end{matrix}\rg),\quad
r_2=\lf(\begin{matrix} \frac{a_{11}}{a_{12}}\vspace{0.4cm}\\
 \frac{a_{21}}{a_{12}}
\end{matrix}\rg)\equiv\lf(\begin{matrix} \frac{\theta_+p_\theta^+}{e_\theta^+}\vspace{0.2cm}\\
\frac{\mu\theta_+u_+^2}{\kappa v_+^2}
\end{matrix}\rg),
\end{equation} where
$p_\theta^+=p_\theta(v_+,\theta_+),\,\,e_\theta^+=e_\theta(v_+,\theta_+)$.
Defining the matrix $B$ by
\begin{equation}\label{2.11.3}
B=\lf(\begin{matrix} 1&\frac{\theta_+p_\theta^+}{e_\theta^+}\\
-\frac{\theta_+p_\theta^+}{e_\theta^+}&\frac{\mu\theta_+u_+^2}{\kappa
v_+^2}
\end{matrix}\rg)\equiv
\lf(\begin{matrix} 1&b_{1}\\
-b_{1}&b_{2}
\end{matrix}\rg),
\end{equation} we have
\begin{equation}\label{2.12.3}
 B^{-1}J_+B=\lf(\begin{matrix} 0 & 0\\  0& \lambda_2
\end{matrix}\rg)=:\Lambda,\quad \quad B^{-1}=\f{1}{b_2+b^2_1}
\lf(\begin{matrix} b_{2}&-b_{1}\\
b_{1}&1
\end{matrix}\rg).
\end{equation}
Setting $Y\equiv(y_1,y_2)^*:=B^{-1}(W-W^+)^*$, we rewrite
\eqref{2.5-3} as
\begin{equation}\label{2.13.3}\begin{aligned}
&Y_x=\Lambda Y+B^{-1}(G(W)-J_+(W-W^+))=:\Lambda Y+H(Y),\quad x>0,\\
&Y(0)=Y_0\equiv B^{-1}(W^--W^+)^*, \quad
\lim_{x\rightarrow\infty}Y(x)=0,
\end{aligned}\end{equation} where
\begin{equation}\label{2.14.3}
H(Y)=B^{-1}\lf(\begin{matrix} \frac{v_+}{\mu
u_+}(p(\hat{v},\hat{\theta})-p^+ -p_v
^+(\hat{v}-v_+)-p_\theta^+(\hat{\theta}-\theta_+))\vspace{0.3cm}\\
\frac{u_+}{\kappa v_+}(e(\hat{v},\hat{\theta})-e^+ -e_v
^+(\hat{v}-v_+)-e_\theta^+(\hat{\theta}-\theta_+))-\frac{u_+^3}{2\kappa
v_+^3}(\hat{v}-v_+ )^2
\end{matrix}\rg).\end{equation}
Let's denote $(\tilde{f}_1,\tilde{f}_2)= BH(Y)$. Using
\eqref{2.14.3} and
\begin{equation}\label{hhhh}
\hat{v}-v_+=y_1+b_1y_2,\,\,\,\hat{\theta}-\theta_+=-b_1y_1+b_2y_2,\end{equation}
we have
\begin{equation}\label{2.15.3}
\begin{aligned}\tilde{f}_1(y_1,y_2)=&\frac{v_+}{\mu
u_+}(p_{vv}^+(\hat{v}-v_+)^2+2p_{v\theta}^+(\hat{v}-v_+)(\hat{\theta}-\theta_+)
+p_{\theta\theta}^+(\hat{\theta}-\theta_+)^2)\\&
+O(|(\hat{v}-v_+,\hat{\theta}-\theta_+)|^3)\\
=&\frac{v_+}{\mu
u_+}(p_{vv}^+-2b_1p_{v\theta}^++b_1^2p_{\theta\theta}^+)y_1^2+O(|y_1|^3+|y_1||y_2|+|y_2|^2),\end{aligned}
\end{equation}
\begin{equation}\label{2.16.3}
\begin{aligned}\tilde{f}_2(y_1,y_2)=&\frac{u_+}{\kappa v_+}
(e_{vv}^+(\hat{v}-v_+)^2+2e_{v\theta}^+(\hat{v}-v_+)(\hat{\theta}-\theta_+)
+e_{\theta\theta}^+(\hat{\theta}-\theta_+)^2)\\&
-\frac{u_+^3}{2\kappa
v_+^3}(\hat{v}-v_+ )^2+O(|(\hat{v}-v_+,\hat{\theta}-\theta_+)|^3)\\
=&\frac{u_+}{\kappa
v_+}(e_{vv}^+-2b_1e_{v\theta}^++b_1^2e_{\theta\theta}^+-\frac{u_+^2}{2
v_+^2})y_1^2+O(|y_1|^3+|y_1||y_2|+|y_2|^2),\end{aligned}
\end{equation}
where $p_{vv}^+=p_{vv}(v_+,\theta_+)$,
$e_{v\theta}^+=e_{v\theta}(v_+,\theta_+)$ and so on.
By  \eqref{2.12.3}, \eqref{2.13.3}, \eqref{2.15.3} and
\eqref{2.16.3}, we have
\begin{equation}\label{2.16-3}\begin{cases}
y_{1x}=-a_+y_1^2 +f_1(y_1,y_2),\\
y_{2x}=\lambda_2y_2 +f_2(y_1,y_2),\quad Y=(y_1,y_2)^*,\\
Y(0)=Y_0, \quad \lim_{x\rightarrow\infty}Y(x)=0,
\end{cases}\end{equation}
where
\begin{equation}\label{2.17-3}\begin{aligned}
&\begin{aligned}a_+&=-\frac{v_+b_2}{\mu
u_+(b_2+b_1^2)}\lf(p_{vv}^+-2b_1p_{v\theta}^++b_1^2p_{\theta\theta}^+\rg)\\
&+\frac{u_+b_1}{\kappa
v_+(b_2+b_1^2)}\lf(e_{vv}^+-2b_1e_{v\theta}^++b_1^2e_{\theta\theta}^+-\frac{u_+^2}{2
v_+^2}\rg),\end{aligned}\\
&f_1(y_1,y_2)=(b_2+b_1^2)^{-1}(b_2\tilde{f}_1-b_1\tilde{f}_2)+a_+y_1^2=O(|y_1|^3+|y_1||y_2|+|y_2|^2),\\
&f_2(y_1,y_2)=(b_2+b_1^2)^{-1}(b_1\tilde{f}_1+\tilde{f}_2)=O(|y_1|^2+|y_1||y_2|+|y_2|^2).
\end{aligned}\end{equation}
By \eqref{2.11.3}, \eqref{1.8-3},  $u_+<0$,
$e_{v\theta}(v,\theta)=\theta p_{\theta\theta}(v,\theta)$ and
\eqref{1.15-3}, we have $a_+>0$. Therefore, applying the center
manifold theory(see \cite{C 81}), by the same lines as in subsection
2.1 of \cite{KNNZ 10}, there exist a local center manifold
$y_2=h^c(y_1)$ and a local stable manifold $y_1=h^s(y_2)$ such that
if the data $(y_{01},y_{02})$ satisfies
$$y_{02}=h^c(y_{01})\quad \text{and}\quad y_{01}\geq h^s(y_{02}),$$
then the problem \eqref{2.16-3} has a unique smooth solution
satisfying
 \begin{equation} \label{2.19-3}
 |\partial_x^kY(x)|\leq\frac{C\sigma^{k+1}}{(1+\sigma x)^{k+1}}+
 C\sigma\exp(-cx),\quad k=0,1,2,
\end{equation}
where $\sigma=|Y_0|$ and $C$ is a positive constant independently of
$x, \sigma$. Moreover, we have
\begin{equation} \label{2.20-3}\begin{aligned}
& y_{1x}=-a_+z^2(x)+O(z^3(x)+\sigma\exp(-cx)),\\
 & y_{2x}=O(z^3(x)+\sigma\exp(-cx)).
\end{aligned}\end{equation}  where $z(x)$ is given as a solution to the equation
$$z_x=-a_+z^2 +f_1(z,h^c(z))$$
and satisfies
\begin{equation} \label{2.22}
0<c\frac{\sigma}{1+\sigma x}\leq z(x)\leq C\frac{\sigma}{1+\sigma
x},\quad |\partial_x^kz(x)|\leq C\frac{\sigma^k}{(1+\sigma
x)^k},\,\,k=1,2.
\end{equation}
Using \eqref{hhhh}, \eqref{2.19-3}, $v_+>0$ and $\theta_+>0$, it is
easy to check that there is a positive constant $\epsilon_0$ such
that if $\sigma=|Y_0|\leq\epsilon_0$, then
 \begin{equation} \label{2.21-3}
 \underline{m}\leq \hat{v}(x), \hat{\theta}(x)\leq\overline{m},
\end{equation}where $\underline{m}$ and $\overline{m}$ are positive constants independently of $x$.
By \eqref{hhhh}, \eqref{2.19-3}, \eqref{2.20-3}, \eqref{2.21-3},
\eqref{2.3-3} and  $u_+<0$, we get \eqref{1.14-3-3} and
\eqref{1.15-3-3}.  The proof of Proposition \ref{theo 1.1-3} is completed.

\

\noindent {\textbf{Acknowledgments:} The authors would like to thank the referees for their valuable comments. The corresponding author Xiaoding Shi and other two authors
 were  partially supported by National Natural Sciences Foundation of China No. 11671027, 11471321 and 11371348.}

 \end{document}